\documentclass[12pt]{amsart}
\usepackage{amsfonts}

\theoremstyle{definition}

\theoremstyle{remark}

\newcommand{\const}{\mathop{\rm const}\limits}

\newcommand{\supp}{\mathop{\rm supp}\limits}

\begin{document}

\begin{center}

{\bf  Composition operators between two different \\
 bilateral grand Lebesgue spaces} \par

\vspace{3mm}
{\bf E. Ostrovsky}\\

e - mail: galo@list.ru \\

\vspace{3mm}

{\bf L. Sirota}\\

e - mail: sirota3@bezeqint.net \\

\vspace{4mm}

Department of Mathematics and Statistics, Bar-Ilan University, 59200, Ramat Gan, Israel. \\

\vspace{5mm}

 Abstract. \\

\vspace{4mm}

{\it In this paper we consider composition operator generated by
nonsingular measurable   transformation
 between two different Grand Lebesgue Spaces  (GLS);
 we investigate the boundedness,
compactness and essential norm of composition operators.}

\vspace{4mm}

\end{center}

\vspace{4mm}

2000 {\it Mathematics Subject Classification.} Primary 37B30,
33K55; Secondary 34A34, 65M20, 42B25.\\

\vspace{4mm}

  {\it Key words and phrases:} norm, Grand and ordinary Lebesgue Spaces, Lorentz and Orlicz spaces,
  rearrangement invariantness, composition linear  operator, H\"older inequality, fundamental  function,
 boundedness, compactness, exact estimations, Radon - Nikodym  derivative. \\

\vspace{3mm}

\section{Introduction}

\vspace{3mm}

  Let $  (X = \{x\}, A, \mu)  $  be measurable space equipped with a non-zero sigma-finite measure $  \mu. $
 Denote by $ M_0  = M_0(X) $ the linear set of all numerical measurable functions  $ f: X \to R. $
  Let also $  \xi = \xi(x)  $ be measurable function from $  X   $ to itself: $ \xi: X \to X.  $ \par

 \vspace{3mm}

 {\bf Definition 1.1.} \par
  The linear operator  $  U_{\xi} [f] = U_{\xi} [f](x) $ defined may be on some Banach subspace $  B_1 $ of the space  $  M_0(X) \ $
to another, in general case, space $  B_2 $   by the formula

 $$
 U_{\xi} [f](x) = f(\xi(x))  \eqno(1.1)
 $$
is said to be {\it  composition operator }  generated by $  \xi(\cdot). $ \par

\vspace{3mm}

 Many  important properties on these operators acting on different spaces  $ B_1 $  with values in $  B_2: $ Lebesgue, Lorentz,
mainly  Orlicz  spaces etc., namely: boundedness, compactness, the exact values of norm  are investigated, e.g. in
 \cite{Arora1} - \cite{Zheng1}. \par

 Applications other than those mentioned appears in ergodic theory, see \cite{Anosov1} - \cite{Petersen1}.\par

 \vspace{3mm}

{\bf  Our purpose in this short article is investigation of these operators in Grand  Lebesgue spaces (GLS).} \par

\vspace{3mm}

 Another operators acting in these spaces: Hardy, Riesz, Fourier, maximal, potential etc.   are investigated, e.g.  in
 \cite{Liflyand1}, \cite{Ostrovsky100}, \cite{Ostrovsky101}, \cite{Ostrovsky3}, \cite{Ostrovsky4}. \\

\bigskip

\section{Grand Lebesgue Spaces (GLS). }

\vspace{3mm}

We recall first of all here  for reader conventions some definitions and facts from
the theory of GLS spaces.\par

Recently, see
\cite{Fiorenza1}, \cite{Fiorenza2},\cite{Ivaniec1}, \cite{Ivaniec2}, \cite{Jawerth1},
\cite{Karadzov1}, \cite{Kozatchenko1}, \cite{Liflyand1}, \cite{Ostrovsky1}, \cite{Ostrovsky2} etc.
 appear the so-called Grand Lebesgue Spaces GLS
 $$
 G(\psi) = G = G(\psi ; A;B);  \ A;B = \const; \ A \ge 1, \ B \le \infty
 $$
spaces consisting on all the measurable functions $ f : X \to R  $ with finite norms

$$
||f||G(\psi) \stackrel{def}{=} \sup_{p \in (A;B)} \left[\frac{|f|_p}{\psi(p)} \right], \eqno(1.2)
$$
where as usually

$$
|f|_p = \left[\int_X |f(x)|^p \ \mu(dx) \right]^{1/p}.
$$

 Here $ \psi = \psi(p), \ p \in (A,B) $ is some continuous positive on the {\it open} interval $ (A;B) $ function such
that

$$
\inf_{p \in(A;B)} \psi(p) > 0. \eqno(1.3)
$$

 We define formally  $ \psi(p) = +\infty, \ p \notin [A,B]. $\par

We will denote
$$
\supp(\psi) \stackrel{def}{=} (A;B).
$$

The set of all such a functions with support $ \supp(\psi) = (A,B) $ will be denoted by  $  \Psi(A,B). $  \par

This spaces are rearrangement invariant; and are used, for example, in
the theory of Probability, theory of Partial Differential Equations,
 Functional Analysis, theory of Fourier series,
 Martingales, Mathematical Statistics, theory of Approximation  etc. \par

 Notice that the classical Lebesgue - Riesz spaces $ L_p $  are extremal case of Grand Lebesgue Spaces;
the exponential Orlicz spaces are the particular cases of Grand Lebesgue Spaces, see
 \cite{Ostrovsky2},  \cite{Ostrovsky100}.  \par

  Let a function $  f:  X \to R  $ be such that

 $$
 \exists (A,B): \ 1 \le A < B \le \infty \ \Rightarrow  \forall p \in (A,B) \ |f|_p < \infty.
 $$
Then the function $  \psi = \psi(p) $ may be naturally defined by the following way:

$$
\psi_f(p) := |f|_p, \ p \in (A,B). \eqno(1.4)
$$

\hfill $\Box$

\bigskip

\bigskip

\section{Main result. }

\vspace{3mm}

{\bf 0.} Recall that the homothety   operator $ T_{\lambda}[g]   $ for arbitrary numerical function $  g $
is defined as ordinary:

$$
T_{\lambda}[g]  = \lambda \cdot  g, \ \lambda = \const.
$$

\vspace{3mm}

 {\bf 1.} We suppose first of all that he source function $  f = f(x) $ belongs to some Grand Lebesgue Space
$   G \psi:   $

$$
|f|_p \le ||f||G\psi \cdot \psi(p), \ p \in (A,B). \eqno(3.0)
$$

 Evidently, the function $ \psi(\cdot) $ may be natural  for  $ f(\cdot), $ if it is non-trivial. In this case
 $ ||f||G\psi  = 1. $\par

\vspace{3mm}

{\bf 2. } Suppose also that the distribution $  \xi $ is absolute continuous relative the measure $  \mu $ with
the Radon - Nikodym  derivative $  h = h(x): $

$$
 \mu \{ y: \xi(y) \in A  \} = \int_A h(x) \ \mu(dx). \eqno(3.1)
$$
 The necessity of this condition is discussed in  \cite{Estaremi1}, \cite{Estaremi2}.\par

  Assume further that the function $ h = h(x) $ belongs to some $ G \theta  $ space:

$$
|h|_q \le ||h||G\theta \cdot \theta(q), \ q \in (a,b). \eqno(3.2)
$$

 In particular, the function $ \theta(\cdot) $ may be natural  for  $ h(\cdot), $ if of course it is non-trivial. In this case
 $ ||h||G\theta  = 1. $\par

\vspace{3mm}

{\bf 3.} Let us introduce a new function $ \nu(p)  = \psi \odot \theta(p) =  [\psi \odot \theta](p) $  by the following way:

$$
 \nu(p) = \psi \odot \theta(p) := \inf_{\alpha > 1} \left\{  \psi(\alpha p) \cdot \left(T_{||h||G\theta}  \right)[ \theta]^{1/p} \left( \frac{\alpha}{\alpha - 1} \right)  \right\}.
  \eqno(3.3)
$$

 {\it It will be presumed that it is finite for some non empty interval } $ p \in (c,d), \ 1 \le c < d \le \infty.  $ \par

\vspace{3mm}

{\bf Theorem 3.1.} {\it  We propose under all formulated in this section conditions: }

$$
||U_{\xi}[f]|| G (\psi \odot \theta) \le ||f||G \psi. \eqno(3.4)
$$

\vspace{3mm}

{\bf Proof.}  We have using H\"older inequality and denoting for brevity \\
 $ \beta = \alpha' = \alpha/(\alpha - 1), \ \alpha > 1: $

$$
| U_{\xi}[f]|_p^p = \int_X |f(\xi(x)) |^p \ \mu(dx) = \int_X |f(x)|^p \ h(x) \ \mu(dx) \le
$$

$$
\left[\int_X |f(x)|^{ \alpha p } \ d \mu \right]^{1/\alpha} \cdot \left[ \int_X |h(x)|^{\beta} \ d \mu \right]^{1/\beta}=
 |f|_{ \alpha p}^p \cdot |h|_{\beta},       \eqno(3.5)
$$
or equally

$$
| U_{\xi}[f]|_p \le |f|_{ \alpha p} \cdot |h|_{\beta}^{1/p}. \eqno(3.6)
$$

 It follows from the direct definition of GLS that

 $$
  |f|_{ \alpha p} \le ||f||G\psi \cdot \psi(\alpha p); \ |h|_{\beta} \le ||h||G\theta \cdot \theta(\beta), \eqno(3.7)
 $$
and we get after substituting into (3.6):

$$
| U_{\xi}[f]|_p \le  ||f||G\psi  \cdot  \left[ \psi(\alpha p) \  [||h||G\theta \cdot \theta(\alpha/(\alpha-1) )]^{1/p} \right].
\eqno(3.8)
$$

 We deduce after minimization over $ \alpha; \ \alpha > 1 $

$$
| U_{\xi}[f]|_p  \le  \nu(p),
$$
which is equivalent to the assertion of theorem 3.1.\par

\vspace{3mm}

{\bf Remark 3.1.}  We conclude in the case of natural pick of both the functions $ \psi(\cdot) $  and $ \theta(\cdot): $

$$
 \nu(p) = \psi \odot \theta(p) := \inf_{\alpha > 1} \left\{  \psi(\alpha p) \cdot [ \theta]^{1/p} \left( \frac{\alpha}{\alpha - 1} \right)  \right\}.
  \eqno(3.9)
$$
and correspondingly

$$
||U_{\xi}[f]||G\nu \le 1. \eqno(3.10)
$$

\bigskip

\section{Examples. }

\vspace{3mm}

{\bf Example 1.} (Exactness.) \par
\vspace{3mm}

 Let us show that the (implicit) coefficient "1" in the right-hand side of inequality (3.4) is in general case non-improvable. \par
   It is  sufficient to consider  a very simple example when $ \xi(x) = x;  $ then $ h(x) = 1. $\par

\vspace{3mm}

{\bf Example 2.} (Power substitution.) \par
\vspace{3mm}

 Let $ X = (0,1) $ with Lebesgue measure $ d \mu = dx.  $   Let now $ f: X \to R  $ be any function from the space $  G\psi. $
 We choose $ \xi(x) = x^m, \ m = \const > 0, $ (not necessary to be integer), so that

 $$
 g(x) = g_m[f](x) \stackrel{def}{=} f(x^m). \eqno(4.1)
 $$

 Define a new function

$$
\psi^{(m)} (p) = \inf_{ \alpha > \min(1,m) }
\left\{ \left[ m^{-1/\alpha} \ (\alpha -1)^{ (\alpha - 1)/\alpha } \ (\alpha - m)^{-1 - 1/\alpha} \right] \cdot \psi(\alpha p) \right\},
\eqno(4.2)
 $$
and suppose this function is non-trivial:  $ \exists (c,d): 1 \le c < d \le \infty \ \Rightarrow  \psi^{(m)} (p) < \infty. $\par

 It follows from theorem 3.1 that

 $$
 ||g_m[f]||G\psi^{(m)} \le  1 \times ||f||G\psi, \eqno(4.3)
 $$
and the constant "1" in (4.3) the best possible. \par

\vspace{4mm}
{\bf Example 3.} (Counterexample). \par
\vspace{3mm}

  Let again $ X = (0,1)  $ and  define

 $$
 f(x) = x^{-1/2}, \ \hspace{6mm} \xi(x) = x^3. \eqno(4.4)
 $$
 Then

 $$
 |f|_p  = \left[ \frac{2}{2-p}   \right]^{1/p} =: \psi(p), \ 1 \le p <2; \eqno(4.5)
 $$

$$
h(z) = 3^{-1} z^{-2/3},  \ 0 < z \le 1; \hspace{5mm} |h|_q = 3^{1/p - 1} \cdot (3 - 2p)^{-1/p} =: \theta(p), \ 1 \le q < 3/2; \eqno(4.6)
$$
or equally $ f(\cdot) \in G\psi, \ h \in G\theta,  $ but the superposition function $ g(x) = f(\xi(x)) = x^{-3/2}  $ does not
belong to any $ L_p(X) $ space with $ p \ge 1.  $ \par

\vspace{4mm}
{\bf Example 4.} (Linear substituting). \par
\vspace{3mm}

 Let here $  X = R^d  $ with Lebesgue measure and $ f: R^d \to R  $ be some function  belonging to the space $  G \psi. $
 Let also $  A $ non degenerate linear map (matrix) from $  R^d $ to $  R^d. $ \par

 Define an operator of a view

 $$
 V_A[f] = f(Ax). \eqno(4.7)
 $$
 Obviously,

$$
| V_A[f]|_p^p = \int_{R^d} |f(A x)|^p \ dx = \int_{R^d}  |\det(A)|^{-1} \  |f(y)|^p \   dy =
$$

$$
|\det(A)|^{-1} \ |f|_p^p,
$$
or equally

$$
|V_A[f]|_p = |\det(A)|^{-1/p}  \ |f|_p \eqno(4.8)
$$
and following

$$
|V_A[f]|_p \le |\det(A)|^{-1/p} \cdot ||f||G\psi   \cdot \psi(p). \eqno(4.9)
$$
 Let the function $  \psi(\cdot) $ be factorable:

$$
\psi(p) = \frac{\zeta(p)}{ \tau(p)}, \ p \in (A,B), \eqno(4.10)
$$
where both the functions  $ \zeta(\cdot), \ \tau(\cdot) $ are from the set $  G\Psi, $ i.e. satisfy all the conditions imposed
on the  function $  \psi(\cdot). $ We deduce after dividing the inequality (4.9) on the function $ \zeta(p): $

$$
\frac{|V_A[f]|_p}{\zeta(p)} \le ||f||G\psi \cdot \frac{[\det(A)]^{-1/p}}{\tau(p)}.\eqno(4.11)
$$

 Recall now that the fundamental  function  $ \phi(G\tau, \ \delta), \ 0 \le \delta \le \mu(X)  $
 for the Grand Lebesgue Space   $  G \tau  $ may be calculated by the formula

$$
\phi(G\tau, \ \delta) = \sup_{p \in (A,B)} \left[ \frac{\delta^{1/p}}{\tau(p)}  \right].
$$

 This notion play a very important role in the theory of operators, Fourier analysis etc., see \cite{Bennett1}.
The detail investigation of the fundamental  function  for GLS is done in \cite{Liflyand1},  \cite{Ostrovsky2}. \par

\vspace{3mm}

 Taking the maximum over $  p; \ p \in (A,B) $ from both the sides of inequality (4.11), we  get to the purpose of this
subsection:  under condition (4.10)

 $$
 ||V_A[f]|| G \zeta  \le ||f||G\psi \cdot \phi(G \tau, |\det(A)^{-1}).  \eqno(4.12)
 $$

\bigskip

\section{Compactness of the composition operator on GLS. }

\vspace{3mm}

 We suppose in this section that the measure $  \mu  $ is diffuse; this imply by known definition that for
arbitrary measurable set $ A $  having positive measure:  $  \mu(A) > 0  $ there exists a measurable subset $ B: \ B \subset A  $
for which $ \mu(B) = \mu(A)/2. $\par

\vspace{3mm}

 We investigate in this section  the compactness of operator $ U_{\xi}[\cdot] $  as an operator acting from source Grand Lebesgue
Space $  G\psi $  into {\it another} GLS $ G\gamma.  $\par
 For Orlicz's, Lebesgue-Riesz, Lorentz  and other Banach  spaces this question is considered in
  \cite{Arora1}, \cite{Cui1},  \cite{Estaremi1},  \cite{Estaremi2},
  \cite{Gupta1}, \cite{Kumar1}, \cite{Singh1}, \cite{Takagi1} -  \cite{Takagi4}.

 \vspace{3mm}

 Let the function $ \nu = \psi \odot \theta, \ \supp \nu  = (a,b)  $ be at the same as defined  in (3.3). \par

\vspace{3mm}

{\bf Theorem 5.1.} {\it  Assume that   }

$$
\sup_{ p \in (a,b)} \left[ \frac{\nu(p)}{\gamma(p)} \right] < \infty \eqno(5.1)
$$
{\it and moreover}

$$
\lim_{\gamma(p) \to \infty} \left[ \frac{\nu(p)}{\gamma(p)} \right] = 0. \eqno(5.2)
$$
{\it (This condition may be omitted if the function  $ \gamma  = \gamma(p)  $ is  bounded.)    } \\
{\it Then the operator }  $ U_{\xi}(\cdot) $ {\it is compact operator from the space } $  G\psi $  {\it into the space } $ G \gamma.  $  \par

\vspace{3mm}

{\bf Proof.}  Denote

$$
B = \{  f:  \ ||f||G\psi \le 1 \},  \hspace{4mm}  U_{\xi} (B) = \{  U_{\xi} [f], \ ||f||G\psi \le 1 \}.  \eqno(5.3)
$$
 Here the symbol $ B $ stands for unit ball in the space $  G \psi.  $ \par
It is sufficient to prove that the set $ U_{\xi} (B)  $ is compact set in the space $ G \gamma. $
In turn it is sufficient  to prove in accordance with the articles \cite{Ostrovsky6}, \cite{Ostrovsky7} that

$$
\lim_{\gamma(p) \to \infty} \sup_{f \in B} \left[ \frac{|U_{\xi} [f]|_p}{\gamma(p)} \right] = 0.  \eqno(5.4)
$$

 It follows from  theorem 3.1 that the set $  U_{\xi} (B)  $  is also bounded in the $ G \nu $ norm:

 $$
 \sup_{f \in B} || U_{\xi} (f) || G \nu \le 1,
 $$
  or equally

 $$
  \sup_{f \in B} | U_{\xi} (f) |_p  \le \nu(p),
 $$
 therefore the left-hand side of inequality (5.4) may be estimated as follows:

$$
\lim_{\gamma(p) \to \infty} \sup_{f \in B} \left[ \frac{|U_{\xi} [f]|_p}{\gamma(p)} \right] \le
\lim_{\gamma(p) \to \infty} \left[ \frac{\nu(p)}{\gamma(p)} \right]  = 0,
$$
Q.E.D. \\

\end{document}